\documentclass[12pt]{article}
\usepackage{amsfonts,amssymb,amsthm,amsmath,amscd,graphicx}

\newcommand{\vv}{\vec{v}}
\newcommand{\vw}{\vec{w}}

\newtheorem{theorem}{Theorem}

\newtheorem{prop}{Proposition}

\newtheorem{defi}{Definition}

\title{ The Full Pythagorean Theorem}
\author{Charles Frohman}
\begin{document}
\maketitle
\begin{abstract} This note motivates  
a version of the generalized pythagorean that says:
if $A$ is an $n\times k$ matrix, then
\[ det(A^tA)=\sum_I det(A_I)^2\] where the sum
is over all $k\times k$ minors of $A$. This is followed
by a proof via an elementary computation in exterior algebra.
The author hopes it is accessible to anybody who has learned
calculus with differential forms.
\end{abstract}

\section{Introduction}
The pythagorean theorem is one of the first theorems of geometry
that people learn.  If a right triangle has legs of length $a$ and $b$ and its hypotenuse has length $c$ then
\[ a^2+b^2=c^2.\]
The Playfair  proof of the Pythagorean theorem is easy to explain, but somehow mysterious.

\vspace{.1in}

\begin{center}
\begin{picture}(128,128)
\includegraphics{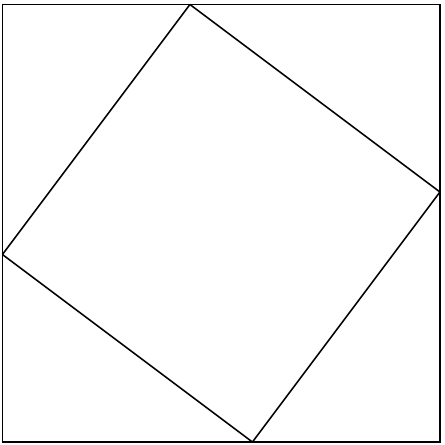}
\put(2,30){$a$}
\put(-25,-10){$b$}
\put(-35,40){$c$}
\put(-100,-10){$a$}
\put(-137,30){$b$}
\put(-100,40){$c$}
\put(-35,85){$c$}
\put(-100,85){$c$}
\put(2,100){$b$}
\put(-137,100){$a$}
\put(-25,132){$a$}
\put(-100,132){$b$}
\end{picture}

\vspace{.25in}
Behold!
\end{center}

The pythagorean theorem leads to the distance formula in
the plane. If $(x_1,y_1), (x_2,y_2)$ are two points in the plane
the distance between them is $\sqrt{(x_1-x_2)^2+(y_1-y_2)^2}$,
which flows from the pythagorean theorem via the diagram below.

\begin{center}
\begin{picture}(128,128)

\includegraphics{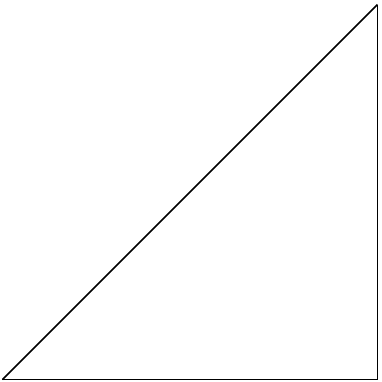}
\put(2,100){$(x_2,y_2)$}
\put(2,0){$(x_2,y_1)$}
\put(-150,0){$(x_1,y_1)$}
\end{picture}
\end{center}

\vspace{.25in}

If you tell someone that the Pythagorean theorem generalizes
to any dimension, often they will think that the generalization is
just the distance formula for points in $\mathbb{R}^n$.  However,
suppose that $T$ is the tetrahedron with vertices at $(0,0,0)$, $(a,0,0)$,
 $(0,b,0)$, and $(0,0,c)$.
 
 \begin{center}
 \begin{picture}(150,173)
 
 \includegraphics{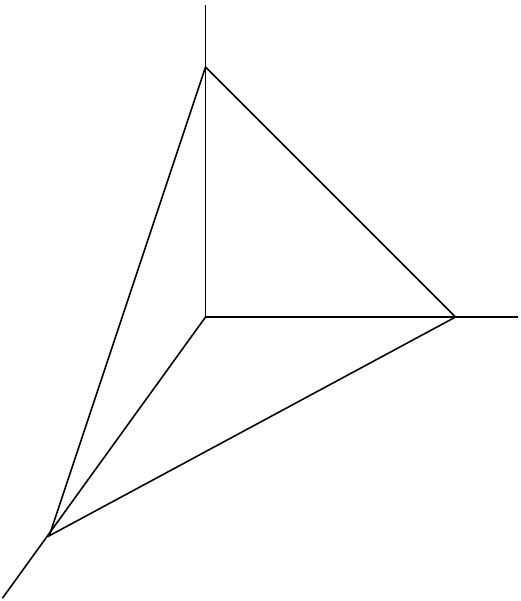}
 \put(-20,70){$(0,b,0)$}
 \put(-90,160){$(0,0,c)$}
\put(-180,20){$(a,0,0)$}

 \end{picture}
 \end{center}
 \vspace{.25in}
 
  You can imagine that the tetrahedron has three legs, which are right triangles that lie in the $xy$-plane, $xz$-plane, and $yz$-plane.  The hypotenuse is the triangle having vertices $(a,0,0)$,
 $(0,b,0)$ and $(0,0,c)$.  The sum of the squares of the areas of the
 three legs is
 \[ \frac{1}{4}a^2b^2+\frac{1}{4}a^2c^2+\frac{1}{4}b^2c^2.\]
 The base of the hypotenuse is $\sqrt{a^2+b^2}$ and its height
 is $\sqrt{\frac{a^2b^2}{a^2+b^2}+c^2}$.  The area of the hypotenuse is one half of its base times its height,
 \[ \frac{1}{2}\sqrt{a^2b^2+a^2c^2+b^2c^2}.\]
 Squaring this gives the sum of the squares of the areas of the legs!  This computation is sometimes referred to as de Gua's theorem after
 Jean Paul de Gua de Malves who was an 18th century French mathematician \cite{We}.

 The statement that this tetrahedron is {\em right}, boils down to the
 legs of the tetrahedron being the orthogonal projections of the hypotenuse into the coordinate hyperplanes.  The shapes we have been using so far, triangles and tetrahedra, are examples of {\em simplices}.  Simplices are well and good, but for the sake of discussing the pythagorean theorem in the language of linear algebra, parallelepipeds are better.  
 
 Suppose that $\vv=(a,b,c)$ and $\vw=(d,e,f)$ are vectors in space.
 The parallelogram spanned by  $\vv$ and $\vw$ is everything of the
 form $s\vv+t\vw$ where $s,t\in [0,1]$.
 The orthogonal projection of the parallelogram into the $xy$-plane
 is spanned by $(a,b)$ and $(d,e)$.  Using the standard area
 formula via determinants its area is $|ae-bd|$. Its orthogonal projection into the $yz$-plane is spanned by $(b,c)$ and $(e,f)$
 and has area $|bf-ce|$.  Finally its projection into the $xz$-plane
 is spanned by $(a,c)$ and $(d,f)$ and has area $|af-cd|$.  The pythagorean theorem says that the square of the area of the parallelogram in space is the sum of the squares of the areas
 of the projections into the coordinate hyperplanes.
 \[ (ae-bd)^2+(bf-ce)^2+(af-cd)^2\]
 but this is just the norm squared of the cross product $(a,b,c)\times (d,e,f)$, which confirms a well known formula for the area of a parallelogram in space.
 
In general, the parallelepiped spanned by vectors $\vv_1,\vv_2,\dots,\vv_k$ is everything of the form $\sum_{i=1}^k\lambda_i\vv_i$ where the $\lambda_i$ vary
over the unit interval $[0,1]$.
 The pythagorean theorem says that if $P$ is a parallelepiped in $\mathbb{R}^n$ spanned by $k$-vectors $\vv_1,\vv_2,\ldots,\vv_k$ then the square of the $k$-dimensional content of $P$, is the sum of the squares of the $k$-dimensional content of the orthogonal projections of $P$ into the $k$-dimensional coordinate hyperplanes in $\mathbb{R}^n$.  It is the point of this note to give a statement in linear algebraic terms of this theorem, and prove it.
 
 \section{Content}
 
 Let $V,W$ be innerproduct spaces, and $L:V\rightarrow W$ a linear map. We can restrict $L$ to get,
 \[ L:\ker{L}^{\perp}\rightarrow \mathrm{im}{L},\]
 where by $\ker{L}^{\perp}$ we mean the subspace of $V$ made up of all vectors that are pependicular to the kernel of $L$ and by $\mathrm{im}{L}$ we mean the image of $L$.  As  $\ker{L}^{\perp}$ and $\mathrm{im}{L}$ are subspaces of innerproduct spaces they are themselves innerproduct spaces. Hence we can choose orthonormal bases for
 them and represent $L$ as  a matrix $l_{ij}$ with respect to those bases.  The matrix $l_{ij}$ is square.
 
 \begin{defi}  The content of $L$, denoted $c(L)$ is the absolute value of the determinant of any matrix representation  $l_{ij}$ 
 of $L$ with respect to orthonormal bases of the perpendicular to its kernel, and its image.
 \end{defi}
 
 The content is not as well behaved as the determinant, for instance it doesn't have a sign.  Also the content of the composition of two linear maps is not in general the product of their contents. However,
 
 \begin{prop} If $L:V\rightarrow W$ and $M:W\rightarrow U$ are linear maps of innerproduct spaces and $\mathrm{im}{L}=\ker{M}^{\perp}$ then $c(M\circ L)=c(M)c(L)$. \end{prop}
 
 \proof Since $\mathrm{im}{L}=\ker{M}^{\perp}$ one of the orthonormal
 bases you use to compute $c(L)$ can be chosen to coincide with one of the orthonormal bases used to compute $c(M)$. If $l_{ij}$ and $m_{jk}$ are the matrices that represent $L$ and $M$ with respect to this choice, then the matrix representing $M\circ L$ is the product of $l_{ij}$ and $m_{jk}$.  The determinant of the product is the product of the determinants. \qed

 If $L:V\rightarrow W$ is a linear map of innerproduct spaces then it has an adjoint $L^*:W\rightarrow V$ defined by requiring $L^*$ to solve the equation,
 \[ <L(\vv),\vw>_W=<\vv,L^*(\vw)>_V\]
 for all $\vv \in V$ and $\vw \in W$ where $<\ ,\ >_V$ and $<\ ,\ >_W$
 denote the innerproducts in $V$ and $W$.  If we have represented $L$ with respect to orthonormal bases of $V$ and $W$ then the matrix of $L^*$ with respect to those bases is just the transpose of the matrix representing $L$. From this it is clear that 
 \[ c(L)=c(L^*).\]
 
 Using Proposition 1, and the fact that the image  of $L^*$ is the perpendicular to the kernel of $L$ we arrive at $c(L\circ L^*)=c(L)^2$,
 and using the fact that $L^{**}=L$ we have $c(L^*\circ L)=c(L)^2$.
 
 \begin{prop} For any linear map of innerproduct spaces $L:V\rightarrow W$, 
 \[ c(L^*\circ L)=c(L\circ L^*)=c(L)^2.\]\end{prop}
 
 Suppose $f:M\rightarrow N$ is an immersion of manifolds and $N$  is Riemannian. Choose local coordinates $x^1,x^2,\ldots,x^k$ at $p\in M$. The matrix
 \[g_{ij}=<df_p(\frac{\partial}{\partial x^i}),df_p(\frac{\partial}{\partial x^j})>_{T_{f(p)}N}\]
 is the first fundamental form of the induced Riemannian metric on $M$ and the associated volume form is
 \[ \sqrt{det(g)} dx^1\wedge \ldots \wedge dx^k.\]
 You can think of $\sqrt{det{g}}(p)$ as the $k$-dimensional content of the parallelepiped
 spanned by the vectors 
 \[df_p(\frac{\partial}{\partial x^1}|_p),\ldots, df_p(\frac{\partial}{\partial x^k}|_p)\]
 Choose an orthonormal basis for
 $T_{f(p)}N$ and let $\vv_i$ be the column vector representing
 \[df_p(\frac{\partial}{\partial x^i}|_p)\] with respect to this basis.  Let $A$ be the $n\times k$-matrix whose columns are the $\vv_i$.  Notice
 that $A^tA=g$, so
 \[ c(A)=\sqrt{det(g)}.\]

In general, you can think of an $n\times k$ matrix $A$ as a list of $k$ column vectors and $c(A)$ as the $k$-dimensional content
of the parallelepiped spanned by those vectors.  Choose $I\subset \{1,2,\ldots,n\}$ with $k$-elements.  We use $|I|$ to denote the cardinality of $I$. By $A_I$ I mean the $k\times k$-matrix made from the rows of $A$ corresponding to the subset $I$.
The parallelepiped spanned by $A_I$ can be thought of as the orthogonal projection of the parallelepiped spanned by $A$ into the $k$-dimensional coordinate hyperplane determined by $I$.  This brings us too:

\begin{theorem}[Pythagorean Theorem] \label{pyth} Let $A$ be an $n\times k$-matrix, then
\[ det(A^tA)=\sum_{I\subset \{1,\ldots, n\} \ |I|=k} det(A_I)^2.\]
That is, the square of the content of the parallelepiped of spanned by $A$ is equal to the sum of the squares of the orthogonal projections of the parallelepiped into the $k$-dimensional coordinate hyperplanes.\end{theorem}

We will prove this theorem after we develop some vocabulary for manipulating determinants.

\section{Exterior Algebra and the Proof}

Let $V$ be a vector space over $k$.  The exterior algebra $\Lambda_*(V)$ of $V$,
is the associative algebra  over $k$ that is the quotient of the free 
unital associative algebra over $k$ on $V$ modulo the relations that $ v\wedge v=0$ for every $v\in V$, where we are denoting the multiplication by a wedge \cite{Sp}.  It is an elementary fact that the relation $v\wedge v=0$ for all $v$ implies that $v\wedge w=-w\wedge v$ for all $v$ and $w$.

We will be working with innerproduct spaces.
Suppose that $e_i$ $i\in \{1,\ldots,n\}$ is an orthogonal basis of $V$. If $I=\{i_1,i_2,\ldots,i_k\}\subset \{1,2,\ldots,n\}$ with $i_1<i_2\ldots <i_k$ let
\[ e_I=e_{i_1}\wedge e_{i_2} \ldots \wedge e_{i_k}.\]
If we let $e_{\emptyset}=1$, then the $e_I$ where $I$ ranges over all subsets of $\{1,2,\ldots, n\}$ is a basis for $\Lambda_*(V)$  In fact,
$\Lambda_*(V)$ is graded by $|I|$.  We let $\Lambda_i(V)$ be the subspace spanned by all $e_I$ with $|I|=i$. Notice $\Lambda_i(V)$ has dimension $\binom{n}{i}$.  We make $\Lambda_*(V)$ into an innerproduct space,  by declaring that the $e_I$ are an orthonormal basis.

If $W$ is any innerproduct space we denote the identity map from $W$ to itself by $Id_W$.
If $L:V\rightarrow W$ is a linear map of innerproduct spaces then
for all $i$, it induces $\Lambda_i(L):\Lambda_i(V)\rightarrow \Lambda_i(W)$, for all $i$.   It is defined by letting
\[ L(\vv_1\wedge\vv_2\ldots \wedge \vv_i)=L(\vv_1)\wedge L(\vv_2)\wedge \ldots \wedge L(\vv_i).\]

The assignment $L\rightarrow \Lambda_i(L)$, is functorial as 
\[ \Lambda_i(Id_V)=Id_{\Lambda_i(V)},\]
and
\[ \Lambda_i(L\circ M)=\Lambda_i(L)\circ\Lambda_i(M).\]
where $L:V\rightarrow W$ and $M:W\rightarrow U$ are any linear mappings
between innerproduct spaces.

\begin{prop}
If we choose orthonormal bases $e_j$ for $V$ and $f_i$ for $W$ and let $l_{ij}$ be the matrix of $f$ with respect to these bases, if $|J|=i$ then 
\[ L(e_J)=\sum_{|I|=i} M_{IJ}f_I\]
where $M_{IJ}$ is the determinant of the $i\times i$ submatrix of $l_{ij}$ whose rows and columns come respectively from the sets $I$ and $J$ 
in order.
\end{prop}

\proof Notice that $L(e_j)=\sum_i l_{ij}f_i$.  Expanding,
\[\Lambda(L)_i(e_J)=(L(e_{j_1})\wedge\ldots \wedge L(e_{j_i}))=\]
\[ (\sum_k l_{kj_1}f_k)\wedge (\sum_k l_{kj_2}f_k)\ldots \wedge  (\sum_k l_{kj_k}f_k)\]
After a lot of cancellation and reordering we arrive at:
\[=\sum_{I\subset \{1,\ldots,n\}}\sum_{\sigma\in S_i} sgn(\sigma)\left(\prod_{m=1}^kl_{i_m,j_{\sigma(m})}\right)f_I,\]
where $S_i$ denotes the symmetric group on $i$ letters, and $sgn(\sigma)$ is the sign of the permutation sigma.  Notice that the sum over $S_i$ of the signed products is a classical formula for the determinant, so it reduces to
\[=\sum_{|I|=i} M_{IJ}f_I\] where the $M_{IJ}$ is the determinant of the  minor whose entries are indexed by
$ij$ with $i\in I$ and $j\in J$.. \qed
\vspace{.1in}

Now suppose that  $dim(V)=dim(W)=k$.  Notice $dim\Lambda_k(V)=dim\Lambda_k(W)=1$, and 
\[ \Lambda_k(L)(e_{\{1,\ldots,k\}})=det(l_{ij})f_{\{1,\ldots,k\}}.\]

\begin{prop} If $L:V\rightarrow W$ is a linear map of innerproduct spaces then the map $\Lambda_i(L^*)$ is equal to $\Lambda_i(L)^*$.
\end{prop}

\proof  Choose orthonormal bases $e_j$ and $f_i$ for $V$ and $W$ respectively.  That way, if $l_{ij}$ is the matrix of $L$ with respect to those bases then the matrix of $L^*$ is $l_{ji}$.  Notice $e_J$ and $f_I$ are orthonormal bases of $\Lambda_i(V)$ and $\Lambda_i(W)$ where we let $I$ and $J$ range over subsets with $i$ elements in their respective index sets. Hence if $\mathcal{L}_{IJ}$ is the matrix of
$\Lambda_i(L)$ with respect to the bases $e_J$ and $f_I$ then the matrix of $\Lambda_i(L^*)$ is $\mathcal{L}_{JI}$.  Here I am using the fact that the determinant of the transpose of a matrix is equal to the determinant of the matrix in computing the coefficients of the matrix for $\Lambda_i(L^*)$. \qed

\vspace{.1in}

We are ready.

\proof {\bf Theorem \ref{pyth} }  Let $A$ be an $n\times k$ matrix,
which we can view as  the matrix of a linear map $\mathcal{A}:\mathbb{R}^k\rightarrow \mathbb{R}^n$ with respect to the standard orthonormal bases.  The determinant of $A^tA$ is expressed as
\[ \Lambda_k(\mathcal{A^*}\circ \mathcal{A})(e_{\{1,\ldots,k\}})=det(A^tA)e_{\{1,\ldots,k\}}.\]  By the functoriality of $\Lambda_k$ the left hand side is
\[ \Lambda_k(\mathcal{A^*})\circ\Lambda_k(\mathcal{A})(e_{\{1,\ldots,k\}}).\]
However, we have
\[\Lambda_k(\mathcal{A})(e_{\{1,\ldots,k\}})=\sum_{I\subset \{1,\ldots n\}\ |I|=k\}} det(A_{I})e_I,\]
where $det(A_I)$ is the determinant of the $k\times k$-matrix whose rows correspond to the rows of $A$ whose index is in $I$.  
Recall,  $\Lambda_k(A^*)=\Lambda_k(A)^*$. Going back to
the definition of adjoint we arrive at
\[ \Lambda_k(A^*)(e_I)=det(A_I)e_{\{1,\ldots,k\}}.\]
Therefore,
\[\Lambda_k(\mathcal{A^*})\circ\Lambda_k(\mathcal{A})(e_{\{1,\ldots,k\}})=\Lambda_k(A^*)(\sum_{I\subset \{1,\ldots n\}\ |I|=k\}}det(A_I)e_{I})=\]
\[\sum_{I\subset \{1,\ldots n\}\ |I|=k\}}det(A_I)\Lambda_k(A^*)(e_{I})=\]
\[ \left(\sum_{I\subset \{1,\ldots n\}\ |I|=k\}}det(A_I)^2\right)e_{\{1,\ldots,k\}}.\]
Putting it all together we have,
\[ det(A^tA)=\sum_{I\subset \{1,\ldots n\}\ |I|=k\}}det(A_I)^2.\] \qed

\section{Epilogue}
This is a  theorem that gets rediscovered over and over again.  Via a web search I found a note by Alvarez \cite{Al} that proves the theorem for right $n$-simplices, but better than that has a nice bibliography with some references to proofs and historical texts.  He cites evidence that the theorem first appeared in a book on analytic geometry by Monge and Hatchet written in the 19th century. 
There is also a paper by Atzema \cite {At} that proves the same theorem I proved here, via a more general result about the determinant of a product of matrices due to Cauchy. I would be surprised if a proof along the lines that
I gave here, didn't appear elsewhere. 

I think the theorem might be of pedagogical interest as it
gives a unified paradigm for integrals to compute arc length
and area that could leave the students in a position to set up
integrals to compute higher dimensional content.


\begin{thebibliography}{0000}
\bibitem{Sp}[Sp] Spivak, Michael  {\em Calculus on Manifolds}, Perseus Books , Cambridge MA, 1965
\bibitem{We}[We] Weisstein, Eric W. {\em de Gua's Theorem.}, From MathWorld--A Wolfram Web Resource, \\ http://mathworld.wolfram.com/deGuasTheorem.html
\bibitem{Al} Alvarez, Sergio {\em Note on an $n$-dimensional Pythagorean theorem} \\ http://www.cs.bc.edu/~alvarez/NDPyt.pdf
\bibitem{At} Atzema, Eisso, {\em Beyond Monge's Theorem}, Mathematics Magazine,
Vol. 73, No. 4 (oct., 2000) pp. 293-296
\end{thebibliography}
\end{document}